\documentclass[12pt]{article}
\usepackage{mathrsfs}
\usepackage{stmaryrd}
\usepackage{amsfonts,amsmath,amssymb,amscd}
\usepackage[vcentermath,enableskew]{youngtab}
\usepackage{shadow}
\usepackage{graphicx}
\usepackage{color}
\usepackage{longtable}
\usepackage{pstricks,multido}
\usepackage{hyperref}
\usepackage{qtree}
\allowdisplaybreaks

\newtheorem{theorem}{Theorem}[section]

\newtheorem{proposition}[theorem]{Proposition}

\newtheorem{conjecture}[theorem]{Conjecture}
\def\qed{\hfill \rule{4pt}{7pt}}

\def\pf{\noindent {\it{Proof.} \hskip 2pt}}

\numberwithin{equation}{section}
\numberwithin{figure}{section}

\begin{document}

\begin{center}
{\large {\bf Vacillating Hecke Tableaux and Linked Partitions}}

\vskip 4mm {\small William Y.C. Chen\\[2pt]
Center for Applied Mathematics\\
Tianjin University\\
Tianjin 300072, P.R. China\\[2pt]
chenyc@tju.edu.cn}

\vskip 3mm {\small Peter L. Guo \\[2pt]
Center for Combinatorics, LPMC-TJKLC\\
Nankai University\\
 Tianjin 300071,
P.R. China\\[2pt]
lguo@nankai.edu.cn  }

\vskip 3mm {\small Sabrina X.M. Pang\\[2pt]
College of Mathematics and Statistics\\
Hebei University of Economics and Business\\
Shijiazhuang, Hebei 050061, P.R. China\\[2pt]
stpangxingmei@heuet.edu.cn  }

\end{center}

\begin{abstract}
We introduce the structure of vacillating Hecke
tableaux, and establish a one-to-one correspondence between
vacillating Hecke tableaux and linked partitions by using the Hecke insertion  algorithm
developed by Buch, Kresch, Shimozono,
Tamvakis and Yong. Linked partitions arise in free probability theory.
Motivated by the Hecke  insertion
 algorithm, we define a Hecke diagram as
a Young diagram possibly with a marked corner.
A vacillating Hecke tableau is defined as a sequence of Hecke diagrams
subject to certain addition and deletion of rook strips. The notion of a rook strip
was  introduced by Buch in the study of the
Littlewood-Richardson rule for stable Grothendieck
polynomials. A rook strip
is a skew Young diagram
with at most one square in each row and column.
We show that the crossing number and the nesting number of a linked partition can be determined by the maximal number of rows and the maximal number of columns of
 the diagrams in the corresponding vacillating Hecke tableau. The proof relies  on a theorem due to
Thomas and Yong  concerning the lengths of the longest strictly increasing and the longest strictly decreasing subsequences
in a word. This implies that the crossing number and the nesting number have a symmetric joint distribution over linked partitions, confirming a conjecture of  de Mier.
We also prove a conjecture of Kim which states  that  the crossing number and the nesting number  have a symmetric joint distribution over the front representations of
partitions.
\end{abstract}

\noindent
{\bf Keywords:} Linked partition, crossing, nesting, the Hecke insertion algorithm

\vspace{6pt}
\noindent
{\bf AMS Classifications:} 05A05, 05A17, 05C30

\section{Introduction}

The Hecke insertion  algorithm  was  developed  by Buch,   Kresch,  Shimozono,  Tamvakis and   Yong \cite{Buch} in order to expand a stable Grothendieck polynomial in terms of stable Grothendieck polynomials indexed by integer partitions. Stable Grothendieck polynomials were defined  by Fomin and Kirillov  \cite{Fomin} as a limit of ordinary Grothendieck polynomials of Lascoux and Sch\"{u}tzenberger \cite{Lascoux}.
For a permutation $\pi$, let $G_\pi=G_\pi(x_1,x_2,\ldots)$ be the stable Grothendieck polynomial indexed by $\pi$. Fomin  and Kirillov \cite{Fomin} showed that $G_\pi$ can be explained as a weighted counting of  compatible pairs. On the other hand, Buch \cite{Buch-02} showed  that
\begin{equation}\label{equ-1}
G_\pi=\sum_{\lambda} c_{\pi,\lambda} G_{\pi_\lambda},
\end{equation}
where the sum ranges over   integer partitions, and for a partition $\lambda$,  $c_{\pi,\lambda}$ is an integer     and $\pi_\lambda$ is a permutation determined by $\lambda$. Buch \cite{Buch-02}  proved that the polynomial  $G_{\pi_\lambda}$ can be  interpreted combinatorially in terms of set-valued tableaux of shape $\lambda$.
A set-valued tableau of shape $\lambda$ is an assignment of nonempty sets of positive integers to the squares of the Young diagram of $\lambda$ such that the sets
are weakly  increasing along each row and strictly increasing along each column,
where for two sets $A$ and $B$ of
integers, the relation  $A\leq B$ (resp., $A<B$) means that $\max(A)\leq \min(B)$ (resp.,
$\max(A)<\min(B)$).

By developing  the Hecke insertion algorithm, Buch et al. \cite{Buch} constructed
a bijection between compatible pairs and pairs $(T,U)$, where $T$ is an increasing tableau and $U$ is a set-valued tableau of the same shape as $T$. An increasing tableau of shape $\lambda$ is an assignment of positive integers to
 the squares of $\lambda$  such that the numbers are strictly increasing in each row and each column. Using this bijection, Buch et al. \cite{Buch} showed that
 the coefficient $c_{\pi,\lambda}$ in   \eqref{equ-1} equals up to a sign
the number of increasing tableaux of shape $\lambda$ satisfying certain conditions.
The Hecke   algorithm can be viewed as an extension of the Robinson-Schensted   algorithm \cite{Robinson,Schensted} and the  Edelman-Greene  algorithm \cite{E-G}.

In this paper, we introduce the structure of vacillating Hecke
tableaux, and  establish a one-to-one correspondence between
vacillating Hecke tableaux and linked partitions by using the Hecke algorithm.
Motivated by  the Hecke algorithm, we define
 a Hecke diagram as
a Young diagram possibly with a marked corner.
A vacillating Hecke tableau is  a sequence of Hecke diagrams
subject to certain addition and deletion of rook strips.
The notion of a rook strip
was  introduced  by Buch \cite{Buch-02} in the study of the
Littlewood-Richardson rule for stable Grothendieck
polynomials. More precisely, a rook strip
is a skew Young diagram
with at most one square in each row and column.
When the Hecke diagrams are restricted to   Young diagrams and the
rook strips are restricted to single squares,
a vacillating Hecke tableau specializes to an ordinary vacillating
tableaux due to Chen, Deng, Du, Stanley and Yan \cite{Chen-1}.

We show that the crossing number   and the nesting number of a linked partition can be determined  by the maximal number of rows and the maximal number of columns of diagrams  in the corresponding vacillating Hecke tableau. The proof relies on   a theorem due to
Thomas and Yong \cite{Thomas} which states  that
the insertion tableau of a word generated by the Hecke algorithm
determines the lengths of the longest strictly increasing and strictly
decreasing subsequences in the word.

As a consequence, we show that  the crossing number and the nesting number have a symmetric joint distribution over linked partitions.
This confirms a conjecture posed   by de Mier \cite{ANNA2006}. We also show that
the crossing number and the nesting number have a symmetric joint distribution over the
front representations of set partitions. This proves a  conjecture of
Kim \cite{Kim}.

The notion of a linked partition
 was introduced by Dykema \cite{DYKE2005} in the study of  unsymmetrized
T-transforms in free probability theory.  Let $[n]=\{1, 2, \ldots, n\}$.
  Dykema \cite{DYKE2005}   showed that   noncrossing linked partitions of $[n+1]$ are counted by the $n$-th large Schr\"oder number. Chen, Wu and Yan
\cite{Chen-2} found a combinatorial interpretation of this fact by establishing a bijection between
noncrossing linked partitions of $[n+1]$ and   Schr\"oder paths of length $n$.

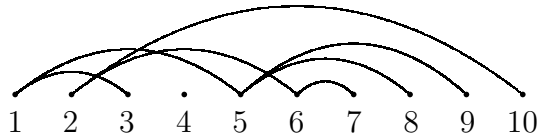
\begin{figure}[ht]
\begin{center}
\setlength{\unitlength}{0.5mm}
\begin{picture}(144,35)

\multiput(2,8)(15,0){10}{\circle*{1.5}}
\put(0,-2){$1$}     \put(12,-2){ $2$}
\put(27,-2){ $3$}    \put(42,-2){ $4$}
\put(57,-2){ $5$}    \put(72,-2){ $6$}
\put(90,-2){$7$}    \put(102.5,-2){ $8$}
\put(117.5,-2){ $9$}   \put(130,-2){ $10$}
\qbezier[500](2,8)(17,20)(32,8) \qbezier[500](2,8)(32,32)(62,8)
\qbezier[500](62,8)(84.5,27)(107,8)\qbezier[500](62,8)(92,35)(122,8)
\qbezier[500](17,8)(47,32)(77,8) \qbezier[500](17,8)(77,55)(137,8)
\qbezier[500](77,8)(84.5,15)(92,8)
\end{picture}
\end{center}
\caption{The linear representation of a linked partition.}
\label{linked-partition}
\end{figure}

A linked partition of $[n]$ is a collection of nonempty subsets $B_1,B_2,\ldots,B_k$ of $[n]$, called blocks,
such that the union of $B_1,B_2,\ldots,B_k$ is $[n]$ and any two distinct blocks are nearly disjoint.
Two distinct blocks $B_i$ and $B_j$ are said to be nearly disjoint if for any $t\in B_i\cap B_j$, one of
the following conditions holds:
\begin{itemize}
\item[(1)] $t=\min(B_i)$, $|B_i|>1$, and $t\neq \min(B_j)$,
\item[(2)] $t=\min(B_j)$, $|B_j|>1$, and $t\neq \min(B_i)$.
\end{itemize}
The linear representation of a linked partition was defined  by Chen, Wu and Yan
\cite{Chen-2}.  For a linked partition $P$ of $[n]$, list the $n$ vertices $1,2,\ldots,n$ in increasing order  on a horizontal line. For a block $B_i=\{a_1,a_2,\ldots, a_m\}$ with
$m\geq 2$ and $a_1<a_2<\cdots<a_m$, draw an arc from $a_1$ to $a_j$ for $j=2,\ldots,m$.
For example, the linear representation of the linked partition $\{\{1,3,5\},\{2,6,10\},\{4\},\{5,8,9\},\{6,7\}\}$ is illustrated in Figure \ref{linked-partition}.
By definition, it is easily checked that  the linear representation of a linked partition of $[n]$ is a simple graph on $[n]$
 such that
for each vertex $i$ there is at most one vertex $j$ with $1\leq j<i$  that is
connected   to $i$, and vice versa.

The crossing number and the nesting number of a linked partition $P$ are defined based on   $k$-crossings and $k$-nestings in the linear representation of $P$, where $k$ is a positive integer.
We use a pair $(i,j)$ with $i<j$ to denote an arc in the linear representation  of  $P$, and we call $i$ and $j$ the left-hand endpoint and the right-hand endpoint of   $(i,j)$, respectively.
We say that $k$ arcs $(i_1,j_1)$,$(i_2,j_2)$,$\ldots$,$(i_k,j_k)$ of $P$ form a $k$-crossing   if
\[i_1<i_2<\cdots<i_k<j_1<j_2<\cdots<j_k,\]
and form   a $k$-nesting   if
\[i_1<i_2<\cdots<i_k<j_k<\cdots<j_2<j_1.\]
The crossing number $\mathrm{cr}(P)$  of $P$ is defined as the maximal number $k$ such that $P$ has a $k$-crossing. Similarly,  the nesting number $\mathrm{ne}(P)$ of $P$ is the maximal number $k$ such that $P$ has a $k$-nesting. For example, for the linked partition in Figure \ref{linked-partition}, we have $\mathrm{cr}(P)=2$ and
 $\mathrm{ne}(P)=3$.

  de Mier \cite{ANNA2006}
posed the following conjecture and showed that it holds for $i=1$.

\begin{conjecture}[\mdseries{de Mier \cite{ANNA2006}}]\label{C-1}
For any positive integers $i$ and $j$, the number of linked partitions $P$ of $[n]$ with $\mathrm{cr}(P)=i$ and
 $\mathrm{ne}(P)=j$ equals  the number of linked partitions $P$ of $[n]$ with $\mathrm{cr}(P)=j$ and
 $\mathrm{ne}(P)=i$.
 \end{conjecture}

Kim \cite{Kim} posed a conjecture on the joint distribution of
the crossing number and the
nesting number of   front representations of set partitions. A set partition of $[n]$ is a collection of mutually disjoint nonempty subsets whose union is $[n]$. Clearly,
a set partition of $[n]$ is a linked partition of $[n]$ such that any two distinct blocks are disjoint. When  $P$ is a set partition, the linear representation of $P$ is   called
 the front representation of $P$ by Kim \cite{Kim}.

\begin{figure}[ht]
\begin{center}
\setlength{\unitlength}{0.5mm}
\begin{picture}(135,35)
\multiput(2,8)(15,0){10}{\circle*{1.5}}
\put(0,-2){$1$}     \put(12,-2){ $2$}
\put(27,-2){ $3$}    \put(42,-2){ $4$}
\put(57,-2){ $5$}    \put(72,-2){ $6$}
\put(90,-2){$7$}    \put(102.5,-2){ $8$}
\put(117.5,-2){ $9$}   \put(130,-2){ $10$}
\qbezier[500](2,8)(17,20)(32,8) \qbezier[500](2,8)(32,32)(62,8)
\qbezier[500](2,8)(54.5,55)(107,8)
\qbezier[500](17,8)(47,32)(77,8) \qbezier[500](17,8)(69.5,55)(122,8)
\qbezier[500](92,8)(114.5,30)(137,8)
\end{picture}
\end{center}
\caption{The front representation of a set partition.}
\label{front-representation}
\end{figure}
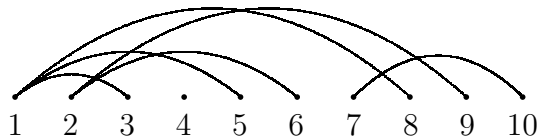

For example, Figure \ref{front-representation} is  the front
representation of
the set partition \[ \{\{1,3,5,8\},\{2,6,9\},\{4\},\{7,10\}\}.\]

\begin{conjecture}[\mdseries{Kim \cite{Kim}}]\label{C-2}
For any positive integers $i$ and $j$, the number of front representations  of set partitions $P$ of  $[n]$ with $\mathrm{cr}(P)=i$ and
 $\mathrm{ne}(P)=j$  equals  the number of front representations of set partitions $P$ of  $[n]$ with $\mathrm{cr}(P)=j$ and
 $\mathrm{ne}(P)=i$.
 \end{conjecture}

We remark  that the crossing number and the nesting number of a set partition have been defined based on the standard representation, see \cite{Chen-1}. The distributions of the crossing number and the nesting number have been extensively studied  in the context  of  fillings of  diagrams, see, for example, \cite{Mier-1,KRAT2006,Yan,Rubey}.

This paper is organized as follows. In Section 2, we present the definition of a vacillating Hecke tableau. In Section 3, we give an overview  of the
Hecke  algorithm developed by Buch et al.\ \cite{Buch}, as well as some properties of this algorithm.
Section 4 provides a bijection between vacillating
Hecke tableaux and linked partitions based on the Hecke algorithm. In particular, we   show that the crossing number and the nesting number of a linked partition are determined  by the maximal number of rows and the maximal number of columns of diagrams in the corresponding vacillating Hecke tableau. As consequences, we confirm Conjecture \ref{C-1} and Conjecture \ref{C-2}.

\section{Vacillating Hecke Tableaux}

In this section, we define
a  vacillating Hecke tableau  as  a sequence of Hecke diagrams
subject to certain addition and deletion of rook strips.
We  first give the notion of a Hecke diagram.
 Let $\lambda=(\lambda_1,\lambda_2,\ldots,\lambda_\ell)$ be a partition of a positive integer $n$,
 that is, $\lambda_1\geq \lambda_2\geq \cdots\geq \lambda_\ell>0$ and
$\lambda_1+\lambda_2+\cdots+\lambda_\ell=n$.
The  Young diagram of $\lambda$ is a left-justified array of squares  with $\lambda_i$ squares in  row $i$.

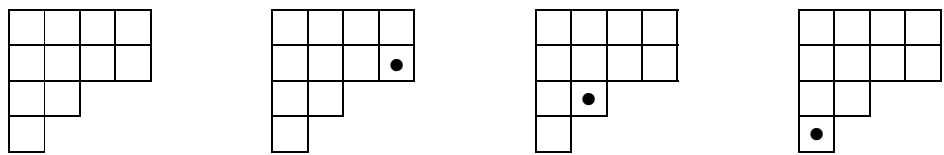
\begin{figure}[h,t]
\setlength{\unitlength}{0.35mm}
\begin{center}
\begin{picture}(160,52)
\put(-100,-8){\yng(4,4,2,1)}
\put(0,-8){\yng(4,4,2,1)}\put(100,-8){\yng(4,4,2,1)}\put(200,-8){\yng(4,4,2,1)}
\put(44.5,22){$\bullet$}\put(117.5,9){$\bullet$}\put(204,-4.5){$\bullet$}
\end{picture}
\end{center}\caption{Hecke diagrams with underlying Young diagram $(4,4,2,1)$.}
\label{Mark-diagram}
\end{figure}

A Hecke diagram is defined as a Young diagram possibly with a marked corner.
For example,  Figure \ref{Mark-diagram} gives  the  Hecke diagrams  whose underlying
Young diagram is $(4,4,2,1)$, where we use a bullet to indicate a marked corner.
We call a Hecke diagram an ordinary diagram if it does
not have a marked corner, and a marked diagram if it has a marked corner.
When $\lambda$ is a marked diagram with a marked corner $c$,  we also write $\lambda$ as a pair   $(\mu,c)$, where $\mu$ is the underlying Young diagram  of $\lambda$.

\begin{figure}[h,t]
\setlength{\unitlength}{0.35mm}
\begin{center}
\begin{picture}(180,52)
\put(120,5){\young(:::\ ,:\ \ ,\ )}\put(0,5){\young(:::\ ,::\ ,\ )}
\put(22,-10){\small{(a)}}\put(142,-10){\small{(b)}}
\end{picture}
\end{center}\caption{Examples of skew diagrams.}
\label{skew-shape}
\end{figure}
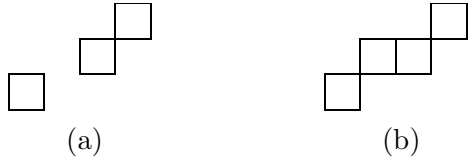

To define
a  vacillating Hecke tableau, we need  the notion of a rook strip    introduced  by Buch \cite{Buch-02}. For two Young diagrams $\lambda$ and $\mu$ such that $\mu$
is contained in $\lambda$, the skew diagram $\lambda/\mu$ is the collection of squares of $\lambda$ that are outside $\mu$.
A rook strip is a skew diagram with at most
one square in each row and  column. For example, in Figure \ref{skew-shape},
the skew diagram (a) is   a rook strip, but (b) is  not a rook strip.

\begin{figure}[h,t]
\setlength{\unitlength}{0.35mm}
\begin{center}
\begin{picture}(420,50)
\put(-5,10){$\emptyset$} \put(15,10){$\emptyset$}
\put(35,10){$\yng(1)$}\put(65,10){$\yng(1)$}\put(95,3){$\yng(2,1)$}
\put(135,3){$\yng(2,1)$}
\put(175,3){$\yng(2,1)$}\put(215,3){$\yng(1,1)$}\put(245,3){$\yng(1,1)$}
\put(275,3){$\yng(1,1)$}\put(305,3){$\yng(1,1)$}\put(335,10){$\yng(1)$}
\put(365,10){$\yng(1)$}\put(395,10){$\emptyset$}\put(415,10){$\emptyset$}

\put(-5,33){0}\put(15,33){1}\put(39,33){2}\put(69,33){3}\put(107,33){4}
\put(147,33){5}\put(185,33){6}\put(218,33){7}\put(250,33){8}
\put(279,33){9}\put(305,33){10}\put(336,33){11}\put(366,33){12}
\put(392,33){13}\put(412,33){14}\put(249,-4){$\bullet$}
\put(99,-4){$\bullet$}\put(38.8,9.7){$\bullet$}
\end{picture}
\end{center}\caption{A vacillating  Hecke tableau of empty shape and length 14.}
\label{va-active-tableaux}
\end{figure}
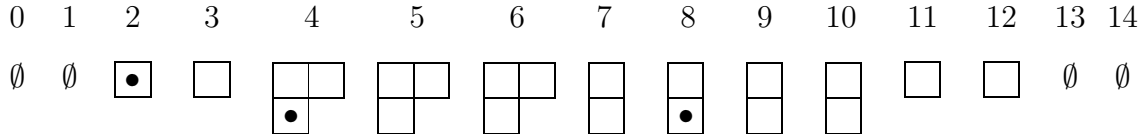

A vacillating Hecke tableau  of
empty shape   and length $2n$ is defined to be a sequence $(\lambda^0,\lambda^1,\ldots,\lambda^{2n})$ of Hecke diagrams   such that
\begin{itemize}
\item[(i)] $\lambda^0=\lambda^{2n}=\emptyset$, and for $1\leq i\leq n$, $\lambda^{2i-1}$ is an  ordinary diagram;
\item[(ii)] If $\lambda^{2i}$ is an ordinary diagram, then $\lambda^{2i-1}$ is an ordinary diagram contained in $\lambda^{2i}$ such that $\lambda^{2i-1}=\lambda^{2i}$ or $\lambda^{2i}/\lambda^{2i-1}$ is a rook strip,  and
    $\lambda^{2i+1}$ is an ordinary diagram contained in $\lambda^{2i}$ such that $\lambda^{2i+1}=\lambda^{2i}$ or $\lambda^{2i}/\lambda^{2i+1}$ is a square;

 If $\lambda^{2i}=(\mu,c)$ is a marked diagram,   then $\lambda^{2i-1}$ is an ordinary diagram   contained in $\mu$ such that $\lambda^{2i-1}=\mu$ or $\mu/\lambda^{2i-1}$ is a rook strip,  and
    $\lambda^{2i+1}=\mu$.

\end{itemize}
As an example,  Figure \ref{va-active-tableaux} illustrates  a vacillating Hecke  tableau of empty shape and length 14.

Notice that when the diagrams of a vacillating Hecke tableau are restricted to
ordinary diagrams and   rook  strips are replaced by
single  squares,  a vacillating Hecke tableau reduces to an ordinary  vacillating tableau  \cite{Chen-1}.

\section{The Hecke insertion algorithm}

In this section, we shall give an overview of the Hecke insertion  algorithm  developed  by Buch et al. \cite{Buch}, as well as a theorem due to Thomas and Yong \cite{Thomas} which
states that the insertion tableau of a word generated by the Hecke insertion algorithm  determines  the
lengths of the longest strictly increasing and strictly  decreasing subsequences in the word.
We also prove a   property concerning the insertion tableau of a word,
which will be used in proof of the main result of this paper, as given in the next section.

The Hecke algorithm  is a procedure to insert  a positive  integer
into  an increasing  tableau, which leads to a representation of a
 word by an increasing tableau and a set-valued tableau.
 Let $\lambda$ be a Young diagram.
An increasing tableau $T$ of shape $\lambda$
is an assignment of  positive integers
to the squares of $\lambda$ such that  the numbers  are strictly increasing in each row and
 column.
Suppose that $U$ is the tableau obtained from $T$ by inserting a positive integer $x$.
Then $U$  is either of the same shape as $T$ or it has
an extra square compared with $T$.
In the case when  $U$ has the same shape  as  $T$, it also contains a special corner
 where the
algorithm terminates and this corner needs to be recorded.
 A parameter $\alpha\in \{0,1\}$ is used to distinguish these two cases.
Thus the  output of the Hecke algorithm when applied to $T$
  is  a triple $(U, c, \alpha)$, where $c$ is a corner of $U$.

The Hecke algorithm can be  described as follows.
Assume that $T$ is an increasing tableau and $x$ is a positive integer.
To insert $x$ into $T$, we begin with the first row of $T$.
Roughly speaking, an element in this row may be
bumped out and then inserted into the next row.
The process is repeated until no more element is bumped out.
More precisely, let $R$ be the first row of $T$. We have the
following two  cases.

\noindent
Case 1: The integer $x$ is larger than or equal to all entries in $R$. If adding $x$ as a new square to the end of $R$ results in
an increasing tableau, then $U$ is the resulting tableau, $c$ is the   corner
where $x$ is added. We set $\alpha=1$ to signify that the corner $c$ is outside the shape of $T$, and the process terminates.
If  adding $x$ as a new square to the end of $R$ does not result in
an increasing tableau, then  let $U=T$, and  $c$ be the corner  at the bottom of the column of $U$ containing the rightmost
square of  $R$. In this case, we set  $\alpha=0$ to indicate that the corner $c$ is inside the shape of $T$, and the process terminates.

\noindent
Case 2: The integer $x$ is strictly smaller than some element in  $R$.
Let $y$ be the leftmost entry in $R$ that is strictly larger than $x$.
If replacing $y$ by  $x$ results in an increasing tableau, then $y$ is bumped out by $x$ and $y$ will be inserted into the next row.
If
replacing $y$ by  $x$ does not result in an increasing tableau, then
 keep the row $R$ unchanged and the element $y$  will also be inserted into the next row.

 We can iterate the above process to insert the element $y$ into the
  next row, still denoted by $R$. Finally, we get the output $(U, c, \alpha)$
  of the insertion algorithm, and we write  $U= (T\xleftarrow{\,\mathrm{H}}x)$
and $(U,c,\alpha)=H(T, x)$.

\begin{figure}[h,t]
\setlength{\unitlength}{0.35mm}
\begin{center}
\begin{picture}(300,145)
\put(24,80){\young(1\ 34,235,45,57)}
\put(0,108){$T\!=$}
\put(80,124){$\leftarrow \!1$} \put(41,123){\bf{2}}
\put(105,108){$\longrightarrow$}

\put(134,80){\young(1234,2\ 5,45,57)}\put(176,110){$\leftarrow \!2$} \put(151,110){\bf{3}} \put(215,108){$\longrightarrow$}

\put(244,80){\young(1234,235,\ 5,57)}\put(274,97){$\leftarrow \!3$} \put(247.5,96){\bf{4}}

\put(24,0){\young(1234,235,35,\ 7)}
\put(-5,28){$\longrightarrow$}
\put(54,4){$\leftarrow \!4$} \put(28,3){\bf{5}}
\put(105,28){$\longrightarrow$}

\put(134,-5){\young(1234,235,35,47,5)}
\put(200,28){$=U$}
\put(145,2){\vector(1,0){20}}\put(168,-1){$c,\ \alpha=1$}

\end{picture}
\end{center}\caption{An example of the Hecke insertion algorithm for $\alpha=1$.}
\label{S-1}
\end{figure}
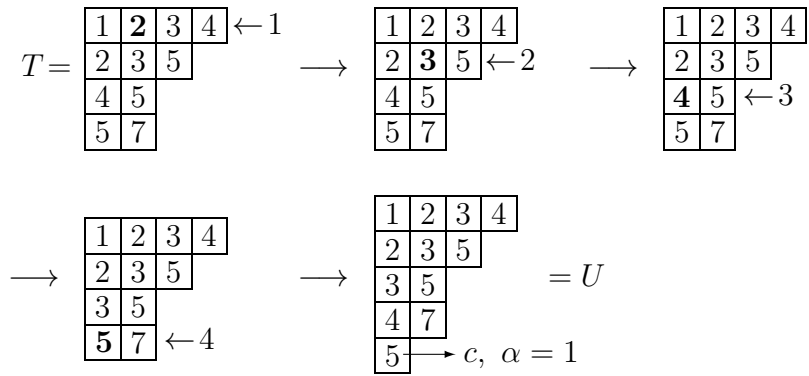

We give two examples to demonstrate the two cases $\alpha=0$ and $\alpha=1$
of the insertion algorithm.
Let $T$ be an increasing  tableau of shape $(4,3,2,2)$ as  given in Figure \ref{IT}.
\begin{figure}[h,t]
\setlength{\unitlength}{0.35mm}
\begin{center}
\begin{picture}(100,55)
\put(40,-8){\young(1234,235,45,57)}
\put(10,20){$T=$}
\end{picture}
\end{center}\caption{An increasing tableau of shape $(4,3,2,2)$.}
\label{IT}
\end{figure}
Let $x=1$. The process to insert $x$ into $T$  is illustrated in Figure \ref{S-1},  where the an element  in boldface represents the entry that is bumped out and  is to be inserted
into the next row.
We see that  the resulting tableau $U$ has one more  square than  $T$, and so
we have $\alpha=1$.

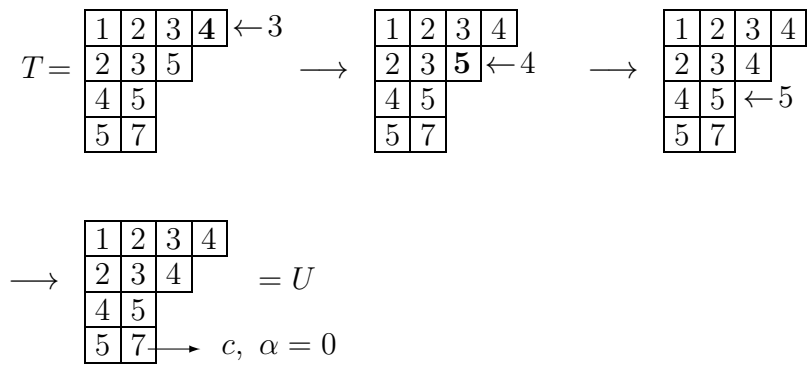
\begin{figure}[h,t]
\setlength{\unitlength}{0.35mm}
\begin{center}
\begin{picture}(300,145)
\put(24,80){\young(123\ ,235,45,57)}
\put(0,108){$T\!=$}
\put(80,124){$\leftarrow \!3$} \put(67,123){\bf{4}}
\put(105,108){$\longrightarrow$}

\put(134,80){\young(1234,23\ ,45,57)}\put(176,110){$\leftarrow \!4$} \put(164,110){\bf{5}} \put(215,108){$\longrightarrow$}

\put(244,80){\young(1234,234,45,57)}\put(274,97){$\leftarrow \!5$}

\put(24,0){\young(1234,234,45,57)}
\put(-5,28){$\longrightarrow$}

\put(90,28){$=U$}
\put(48,5){\vector(1,0){20}}\put(76,3){$c,\ \alpha=0$}
\end{picture}
\end{center}\caption{An example of the Hecke insertion algorithm for $\alpha=0$.}
\label{S-2}
\end{figure}

For $x=3$, we find that the resulting tableau $U$ has the same shape as $T$, and so we have $\alpha=0$, see Figure \ref{S-2}.

The  Hecke algorithm is reversible, see Buch et al. \cite{Buch}. In other words,   give an increasing tableau $U$, a corner $c$   of $U$, and the value of $\alpha$,
there exist a unique increasing tableau $T$ and a unique positive integer $x$  such that $U= (T\xleftarrow{\,\mathrm{H}}x)$.

Thomas and Yong \cite{Thomas} showed that the Hecke algorithm can be used to
determine   the lengths of the longest strictly increasing and strictly decreasing subsequences of a word.
Let $w=w_1w_2\cdots w_n$ be a word of positive integers.
A subword of $w=w_1w_2\cdots w_n$ is a subsequence  $w_{i_1}w_{i_2}\cdots w_{i_k}$, where $1\leq i_1<i_2<\cdots<i_k\leq n$. A subword $w_{i_1}w_{i_2}\cdots w_{i_k}$ is said to be strictly increasing if
$w_{i_1}<w_{i_2}<\cdots< w_{i_k}$, and   strictly decreasing if $w_{i_1}>w_{i_2}>\cdots> w_{i_k}$.
Let $\mathrm{is}(w)$ (resp., $\mathrm{de}(w)$) denote the length
of the longest  strictly   increasing (resp., strictly decreasing)  subwords of $w$.
As shown by Thomas and Yong,  $\mathrm{is}(w)$ and $\mathrm{de}(w)$
are determined by the shape of the insertion tableau of $w$.
The insertion tableau of $w$ is defined by
\[(\cdots((\emptyset\xleftarrow{\,\mathrm{H}} w_1)
\xleftarrow{\,\mathrm{H}} w_2)\xleftarrow{\,\mathrm{H}}\cdots
)\xleftarrow{\,\mathrm{H}} w_n.\]

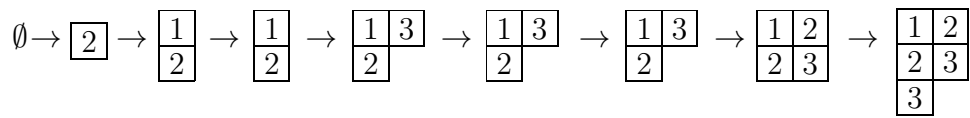
\begin{figure}[h,t]
\setlength{\unitlength}{0.3mm}
\begin{center}
\begin{picture}(450,45)
\put(10,20){$\emptyset$}
\put(18,20){$\rightarrow$}
\put(36,15){\young(2)}\put(56,20){$\rightarrow$}
\put(75,5){\young(1,2)}\put(97,20){$\rightarrow$}
\put(117,5){\young(1,2)}\put(140,20){$\rightarrow$}
\put(161,5){\young(13,2)}
\put(200,20){$\rightarrow$}
\put(220,5){\young(13,2)}
\put(261,20){$\rightarrow$}
\put(282,5){\young(13,2)}
\put(321,20){$\rightarrow$}
\put(340,5){\young(12,23)}
\put(380,20){$\rightarrow$}
\put(402,-10){\young(12,23,3)}
\end{picture}
\end{center}\caption{The insertion tableau of $w=21131321$.}
\label{RSK-bullet}
\end{figure}

For example, let  $w=21131321$. The construction of the insertion tableau
 of $w$ is given in  Figure \ref{RSK-bullet}.

For an increasing tableau $T$, let $\mathrm{c}(T)$ and $\mathrm{r}(T)$ denote the number of columns  and the number of rows  of $T$, respectively.  Using the jeu de taquin algorithm for increasing tableaux developed in \cite{TY}, Thomas and Yong \cite{Thomas} established the following relation.

\begin{theorem} \label{connection}
Let $w$ be a word of positive integers, and $T$ be the insertion tableaux of $w$. Then
$\mathrm{is}(w)=\mathrm{c}(T)$ and $\mathrm{de}(w)=\mathrm{r}(T)$.
\end{theorem}

We observe the following  property of the insertion tableau.

\begin{proposition}\label{deleting}
Let $w=w_1w_2\cdots w_n$ be a word of positive integers, and $k$ be the maximal element appearing in $w$.
Let $w'=a_1a_2\cdots a_m$ be the word obtained from $w$ by deleting the elements equal to $k$.
Assume that $T$ is the insertion tableau of $w$ and $T'$ is the insertion tableau of $w'$.
Then  $T'$ is obtained from $T$ by deleting the squares occupied with $k$.
\end{proposition}

\pf Let $Q$ denote the increasing tableau obtained from $T$ by deleting the squares occupied with the maximal element $k$.
We use induction to prove that $T'=Q$. The claim is obvious when $n=1$. We now assume that
$n>1$ and that the claim holds for $n-1$. Let $P$ be the insertion tableau of $w_1w_2\cdots w_{n-1}$. Here are two cases.

\noindent
Case 1: $w_n=k$.
By the  induction hypothesis, $T'$ is obtained from $P$ by deleting the squares occupied with $k$.
On the other hand, since $w_n=k$ is the maximal element of $w$, we see that  $T=P$ or $T$ is obtained from $P$ by adding a square filled with $k$ at the end of  the first row.
This yields that $T'=Q$.

\noindent
Case 2: $w_n<k$.
 Let $U$ be the insertion
tableau of $a_1\cdots a_{m-1}$. By the induction hypothesis, $U$  is obtained from
$P$ by deleting the squares occupied  with $k$.
In the process of inserting $w_n$ into $P$,
if no entry      equal to $k$ is   bumped out and is  inserted into the next row, then it is clear that $T'=Q$.

Otherwise, there is a  unique entry $k$  in $P$  that is bumped out and is inserted into the next row. Let $c$ be the  square of $P$ occupied with this entry. Note that $c$ is a corner of $P$ since $P$ is increasing and  $k$ is a maximal entry.
Keep in mind that $U$  is obtained from
$P$ by deleting the squares occupied with $k$. Since $T=(P\xleftarrow{\,\mathrm{H}}w_n)$
and $T'=(U\xleftarrow{\,\mathrm{H}}w_n)$, for any square $C$ in $U$,
the entry of $T'$ in $C$ equals the entry of $T$ in $C$.
Consequently, to verify $T'=Q$, it suffices to consider the entry  of $T$ in the corner $c$.
Assume that this entry is equal to $i$.
Here are two subcases.

\noindent
 Case 2.1: $i=k$. In this case, $T'$ has the same shape as $U$. On the other hand, any square of $T$ outside  $U$ is occupied with $k$. So we have $T'=Q$.

\noindent
Case 2.2: $i<k$. In this case, $T'$ has the extra corner $c$ compared with $U$.
 By the construction of the Hecke algorithm, we see that the entry of $T'$ in the corner $c$ also  equals $i$. Notice also  that except for the corner $c$, any square of $T$
outside  $U$ is occupied  with $k$. So we are led to $T'=Q$.
This completes the proof.
\qed

For example, let
$w=21131321$. Then we have $w'=211121$. The insertion tableau $T$ of $w$ is given in Figure \ref{RSK-bullet}.
Meanwhile, the insertion tableau $T'$ of $w'$ is constructed in Figure \ref{deletion}, which
coincides with
the  tableau obtained from  $T$  by deleting the squares occupied  with $3$.
\begin{figure}[h,t]
\setlength{\unitlength}{0.3mm}
\begin{center}
\begin{picture}(315,40)
\put(10,8){$\emptyset$}
\put(18,10){$\rightarrow$}
\put(36,5){\young(2)}\put(56,10){$\rightarrow$}
\put(75,-5){\young(1,2)}\put(95,10){$\rightarrow$}
\put(115,-5){\young(1,2)}\put(140,10){$\rightarrow$}
\put(161,-5){\young(1,2)}
\put(186,10){$\rightarrow$}
\put(206,-5){\young(12,2)}
\put(247,10){$\rightarrow$}
\put(268,-5){\young(12,2)}
\end{picture}
\end{center}\caption{The insertion tableau of $w'=211121$.}
\label{deletion}
\end{figure}

\section{Vacillating Hecke tableaux and linked partitions}

In this section,  we provide a bijection between vacillating
Hecke tableaux and linked partitions. We prove  that the crossing number and the nesting number of a linked partition can be determined  by the maximal number of rows and the maximal number of columns of diagrams in the corresponding vacillating Hecke tableau.
As a consequence, we show that the crossing number and the nesting number have
a symmetric joint distribution over linked partitions with fixed  left-hand endpoints and  right-hand endpoints. This leads to a proof  of  Conjecture \ref{C-1}.
Specializing the bijection to linked partitions containing
 no vertex that is both   a left-hand endpoint and a right-hand endpoint,
 we confirm  Conjecture \ref{C-2}.

To describe our bijection, by  a Hecke tableau we mean an increasing tableau possibly with a marked corner.
In other words, a Hecke tableau is an increasing tableau whose shape is a
Hecke diagram. Let $\lambda$ be a Hecke diagram, and let  $T$  be a Hecke tableau of shape $\lambda$.
When $\lambda=(\mu,c)$ is a marked diagram, we also
express $T$ by a pair $(T',c)$,
where $T'$ is the underlying increasing  tableau of $T$.

Let $V_{2n}$ be the set of  vacillating Hecke tableaux of empty shape and length $2n$.
We now give a description of a bijection $\phi$ from $V_{2n}$ to the set of linked partitions of $[n]$.

Let
$V$ be a vacillating Hecke
  tableau $(\emptyset=\lambda^0,\lambda^1,\ldots,\lambda^{2n}=\emptyset)$ of empty shape and length $2n$.
First, we recursively define a sequence
$(E_0,T_0),(E_1,T_1),\ldots,(E_{2n},T_{2n})$, where for $0\leq i\leq 2n$,
$E_i$ is a set of  edges and $T_i$ is a Hecke tableau of   shape
$\lambda^i$. Let $E_0=\emptyset$, and let  $T_0$ be the
empty tableau. Assume that $i\geq 1$. If $\lambda^i=\lambda^{i-1}$, then we let
$(E_i,T_i)=(E_{i-1},T_{i-1})$. If $\lambda^i\neq\lambda^{i-1}$, $(E_i,T_i)$
 is constructed according to the parity of $i$.

\noindent
Case 1:  $i$ is odd. Let $i=2k-1$. By the definition of a vacillating Hecke tableau, $\lambda^i$ is a  diagram without any marked corner.
Here are two subcases according to whether $\lambda^{i-1}$ is an ordinary diagram.

\noindent
Case 1.1:
$\lambda^{i-1}$ is an ordinary diagram.
 Then $\lambda^i$ is obtained from $\lambda^{i-1}$ by deleting a corner $c$.  Setting $\alpha=1$,
 there are a unique increasing tableau
$T$ and a unique  positive integer $j$ such that $(T_{i-1},c,\alpha)=H(T,j)$.
Let $T_i=T$ and define  $E_i$ to be the set obtained from $E_{i-1}$ by adding the edge $(j,k)$.

\noindent
Case 1.2:
$\lambda^{i-1}=(\mu,c)$ is a marked diagram. So we have $\lambda^i=\mu$.
Setting $\alpha=0$, there are a unique increasing tableau
$T$ and a unique  positive integer $j$ such that $(T_{i-1},c,\alpha)=H(T,j)$.
Let $T_i=T$ and define
$E_i$ to be the set obtained from $E_{i-1}$ by adding the edge $(j,k)$.

\noindent
Case 2:
$i$ is even. Let $i=2k$. We set $E_i=E_{i-1}$. To define $T_i$,
there are two  subcases according to  whether  $\lambda^i$ is an ordinary diagram.

\noindent
Case 2.1:  $\lambda^i$ is an ordinary diagram. Then $\lambda^i/\lambda^{i-1}$ is a rook strip. Define $T_i$ to be the  tableau  obtained
from $T_{i-1}$ by filling the squares
of $\lambda^i/\lambda^{i-1}$ with  $k$.

\noindent
Case 2.2:
$\lambda^i=(\mu,c)$ is a marked diagram.  Then $\mu/\lambda^{i-1}$ is a rook strip.
Let   $T$ be the  tableau of shape $\mu$ that is obtained
from $T_{i-1}$ by filling the squares
of $\mu/\lambda^{i-1}$ with $k$.
Define   $T_i=(T,c)$.

Finally, we define  $\phi(V)$ to be the diagram with $n$ vertices $1,2,\ldots,n$ listed
on a horizontal line such that there is an arc connecting $j$ and $k$ with $j<k$ if and only if
$(j,k)$ is an edge in $E_{2n}$.
By the above construction,  for each vertex $k\in [n]$, there is at most one vertex $j$
with  $j<k$ such that $(j,k)$ is an arc in  $\phi(V)$. Thus $\phi(V)$ is a linked partition of $[n]$.

\begin{figure}[h,t]
\setlength{\unitlength}{0.35mm}
\begin{center}
\begin{picture}(420,100)
\put(-5,60){$\emptyset$} \put(15,60){$\emptyset$}
\put(35,60){$\young(\ )$}\put(65,60){$\young(1)$}\put(95,53){$\young(12,\ )$}
\put(135,53){$\young(12,2)$}
\put(175,53){$\young(12,2)$}\put(215,53){$\young(1,2)$}\put(245,53){$\young(1,\ )$}
\put(275,53){$\young(1,2)$}\put(305,53){$\young(1,2)$}\put(335,60){$\young(2)$}
\put(365,60){$\young(2)$}\put(395,60){$\emptyset$}\put(415,60){$\emptyset$}

\put(-5,80){0}\put(15,80){1}\put(39,80){2}\put(69,80){3}\put(107,80){4}
\put(147,80){5}\put(185,80){6}\put(218,80){7}\put(249,80){8}
\put(279,80){9}\put(305,80){10}\put(336,80){11}\put(365,80){12}
\put(392,80){13}\put(412,80){14}
\put(38.5,59){\bf 1}\put(99,45){\bf 2}
\put(249,45){\bf 2}

\multiput(140,0)(24,0){7}{\circle*{1.5}}
\put(137,-10){\small{1}}     \put(157,-10){ \small{2}}
\put(182,-10){ \small{3}}    \put(205,-10){ \small{4}}
\put(230,-10){ \small{5}}    \put(254,-10){ \small{6}}
\put(281,-10){\small{7}}
\qbezier[500](140,0)(152,8)(164,0)
\qbezier[500](140,0)(164,18)(188,0) \qbezier[500](140,0)(188,38)(236,0)
\qbezier[500](140,0)(200,55)(260,0)
\qbezier[500](164,0)(188,18)(212,0) \qbezier[500](164,0)(224,55)(284,0)

\end{picture}
\end{center}\caption{An illustration of the bijection $\phi$.}
\label{pt}
\end{figure}
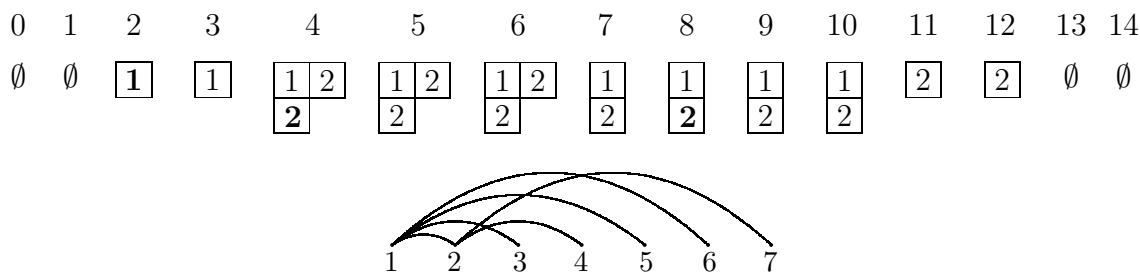

Figure \ref{pt} gives an illustration of the map $\phi$
when applied to the  vacillating Hecke tableau
in Figure \ref{va-active-tableaux},
where an entry  in boldface indicates a  marked  corner of a Hecke tableau.

 It can be checked that
 $\phi$ is reversible, and hence it is a bijection.
Let $P$ be a  linked partition of $[n]$.
To recover the corresponding vacillating Hecke tableau, we  first construct a sequence
$(T_0,T_1,...,T_{2n})$ of Hecke tableaux.
Let $T_{2n}$ be the empty tableau. Suppose that $T_{2i}$ has been constructed, where $1\leq i\leq n$.
We proceed to construct  $T_{2i-1}$ and $T_{2i-2}$. To obtain  $T_{2i-1}$, we have two cases.

 \noindent
Case 1: The shape of $T_{2i}$ is an ordinary diagram. Let $T_{2i-1}$ be the tableau obtained from $T_{2i}$ by deleting the squares (if any) filled with $i$.

\noindent
Case 2: The shape of $T_{2i}$ is a marked diagram. Assume that $T_{2i}=(T,c)$. Let $T_{2i-1}$ be the tableau obtained from $T$ by deleting the squares (if any) filled with $i$.

Now the Hecke tableau $T_{2i-2}$ is obtained from $T_{2i-1}$. If $i$ is not a right-hand endpoint of any arc of $P$,  then  we set $T_{2i-2}=T_{2i-1}$. Otherwise, there is a unique arc $(j,i)$ with $j<i$ of $P$. Assume that $(U,c,\alpha )=H(T_{2i-1},j)$. We set $T_{2i-2}=U$ if $\alpha=1$, and set $T_{2i-2}=(U,c)$ if $\alpha=0$.

Let $\lambda^i$ be the shape of $T_i$. Finally,  the vacillating Hecke tableau $\phi^{-1}(P)$ is given by
\[(\emptyset=\lambda^0,\lambda^1,\ldots,\lambda^{2n}=\emptyset).\]

The following theorem shows that the crossing number and the nesting number
of a linked partition are determined by the diagrams in the corresponding vacillating
Hecke tableau.
For a  vacillating Hecke  tableau $V$, let
$\mathrm{r}(V)$
be  the most number of
rows  in any  diagram $\lambda^i$ of  $V$. Similarly, let $\mathrm{c}(V)$ be
 the most number of
columns  in any  diagram $\lambda^i$ of  $V$.  We have the following relations.

\begin{theorem}\label{key}
Let $V$ be a  vacillating Hecke tableau in $V_{2n}$, and let $P=\phi(V)$.
Then we have $\mathrm{c}(V)=\mathrm{ne}(P)$ and $\mathrm{r}(V)=\mathrm{cr}(P)$.
\end{theorem}

\pf
Let $V=(\lambda^0,\lambda^1,\ldots,\lambda^{2n})$, and let
$(T_0,T_1,\ldots,T_{2n})$ be the  sequence of Hecke tableaux  in
the construction of $\phi$. For $0\leq i\leq 2n$, let $U_i$ be
be the underlying increasing tableau of $T_i$.
We shall  generate a sequence $(w^{(0)},w^{(1)},\ldots,
 w^{(2n)})$ of words such that $U_i$ is the insertion tableau of $w^{(i)}$.

Let $w^{(2n)}$ be the
empty word. Suppose that
 $w^{(i)}$
has been constructed. Then $w^{(i-1)}$ is constructed as follows.
 If $T_{i-1}=T_{i}$, then we let
$w^{(i-1)}=w^{(i)}$. If $T_{i-1}\neq T_{i}$, we have
 two cases.

\noindent
Case 1: $i$ is odd. Let $i=2k-1$.  By the construction of $\phi$,
we see that $U_{i-1}$ is obtained from $U_{i}$
by inserting a unique integer  $j$. Define  $w^{(i-1)}=w^{(i)}\,j$;

\noindent
Case 2: $i$ is even. Let $i=2k$.
Again,   by  the construction of $\phi$, we find that
$U_{i-1}$ is obtained from $U_i$ by deleting the squares (if any)
filled with $k$. Define $w^{(i-1)}$ to be the word obtained from $w^{(i)}$ by removing the elements (if any) equal to $k$.

We proceed by induction to show that  $U_i$ is the insertion tableau of $w^{(i)}$.
 The claim is obvious for $i=2n$. Assume
that the claim is true for  $i$, where $1\leq i\leq 2n$. We wish to prove that
it holds for $i-1$.
If $T_{i-1}=T_i$, then the claim is evident. Let us now consider the case when $T_{i-1}\neq T_i$.  If $w^{(i-1)}$ is generated according to Case 1, then the claim follows directly from the construction of  $\phi$.
If $w^{(i-1)}$ is generated according to Case 2, then the claim is a consequence of Proposition \ref{deleting}. This proves the claim.

Combining the above claim  and Theorem \ref{connection}, we obtain that
\[\mathrm{c}(V)=\max\left\{\mathrm{is}(w^{(i)})\,|\, 0\leq i\leq 2n\right\} \]
and
\[\mathrm{r}(V)=\max\left\{\mathrm{de}(w^{(i)})\,|\, 0\leq i\leq 2n\right\}.\]
It remains to show that  $P$ has
a $k$-crossing (resp., $k$-nesting) if and only if there exists  a word $w^{(i)}$  that contains  a strictly  decreasing (resp., increasing)
subword of length $k$.
We shall only give the proof of the statement concerning the
relationship between  a $k$-crossing and a strictly decreasing subword of length $k$.
The same  argument applies to the $k$-nesting case.
Suppose that $w^{(i)}=a_1 a_2\cdots a_t$  contains a strictly decreasing
subsequence $a_{i_1}\cdots a_{i_k}$  of
 length $k$, where $1\leq i_1<\cdots <i_k\leq t$. By the construction of $\phi$ and the construction
 of the
sequence $(w^{(0)},w^{(1)},\ldots,
w^{(2n)})$,
 we deduce that the vertices $a_{i_1},\ldots,a_{i_k}$ are left-hand endpoints  of $P$.
For $1\leq s\leq k$, let $b_{j_s}$ be the right-hand  endpoint connected  to $a_{i_s}$.
Again, by the construction
 of  $(w^{(0)},w^{(1)},\ldots,
w^{(2n)})$, we see that
\begin{equation}\label{L3}
b_{j_1}>b_{j_2}>\cdots>b_{j_k}.\end{equation}
Hence  the arcs $(a_{i_1}, b_{j_1}),\ldots,(a_{i_k}, b_{j_k})$ form a
 $k$-crossing of $P$.

On the other hand, suppose that $P$ has a $k$-crossing consisting of  arcs
\[(i_1,j_1), (i_2,j_2),\ldots,(i_k,j_k),\] where
$i_1<i_2<\cdots<i_k<j_1<j_2<\cdots<j_k$. By the construction
 of the
sequence $(w^{(0)},w^{(1)},\ldots,
w^{(2n)})$, it is easily checked   that $i_1 i_2 \cdots i_k$ forms a
strictly decreasing subword of $w^{(2j_1-1)}$.
This completes the proof.
 \qed

Conjecture \ref{C-1} and Conjecture \ref{C-2} are consequences of Theorem \ref{key}.
Like the symmetric joint distribution of the crossing number and the
nesting number of ordinary partitions due to Chen, Deng, Du, Stanley and Yan
\cite{Chen-1}, we restrict our attention to fixed sets of the left-hand endpoints   and the right-hand endpoints of a linked partition.
For two subsets $S$ and $T$  of $[n]$, let $L_n(S,T)$ be the set of linked partitions  of $[n]$ such that $S$ is the set of left-hand endpoints and $T$ is the
 set of right-hand endpoints. Note that $L_n(S,T)$ may be empty.

 Let $f_{n,S,T}(i,j)$ be the number of linked partitions $P$ in $L_n(S,T)$ with $\mathrm{cr}(P)=i$ and
 $\mathrm{ne}(P)=j$. We have the following symmetry property.

\begin{theorem}\label{thm-f}
Let $S$ and $T$ be  two subsets of $[n]$.
For any positive integers $i$ and $j$, we have
\[
f_{n,S,T}(i,j)=f_{n,S,T}(j,i).
\]
\end{theorem}

To prove Theorem \ref{thm-f}, we establish an involution on the set $L_n(S,T)$ that exchanges the crossing number and the nesting number of a linked partition. To this end, we  define the conjugate of a Hecke diagram as the transpose of the diagram.

\noindent
\textit{Proof of Theorem \ref{thm-f}.}
Taking the conjugate of every Hecke  diagram  leads to an involution on    vacillating Hecke tableaux of empty shape and length $2n$.
This yields an involution, denoted  $\psi$, on the set  of linked partitions of $[n]$. By Theorem \ref{key}, we find that $\psi$
exchanges the crossing number and the nesting number of a linked partition.
It remains to show that $\psi$ preserves the left-hand endpoints  and the right-hand endpoints of a linked partition.

Let $P$ be a linked partition of $[n]$, and let $(\lambda^0,\lambda^1,\ldots,\lambda^{2n})$ be the corresponding vacillating Hecke tableau.
By the construction of $\phi$, we observe that a vertex $i$ of $P$ is a left-hand endpoint
if and only if $\lambda^{2i}$ has at least one more square than $\lambda^{2i-1}$, and it is a right-hand endpoint   if and only if $\lambda^{2i-2}\neq \lambda^{2i-1}$.
Hence the involution $\psi$  preserves the left-hand   and the right-hand endpoints.
Restricting  $\psi$  to
$L_n(S,T)$ gives an involution on $L_n(S,T)$. This completes the proof. \qed

\begin{figure}[h,t]
\setlength{\unitlength}{0.35mm}
\begin{center}
\begin{picture}(420,40)

\multiput(140,0)(24,0){7}{\circle*{1.5}}
\put(137,-10){\small{1}}     \put(157,-10){ \small{2}}
\put(182,-10){ \small{3}}    \put(205,-10){ \small{4}}
\put(230,-10){ \small{5}}    \put(254,-10){ \small{6}}
\put(281,-10){\small{7}}
\qbezier[500](140,0)(176,25)(212,0) \qbezier[500](140,0)(212,58)(284,0)
\qbezier[500](164,0)(200,25)(236,0)
\qbezier[500](164,0)(176,8)(188,0) \qbezier[500](164,0)(212,40)(260,0)
\qbezier[500](140,0)(152,8)(164,0)
\end{picture}
\end{center}\caption{An example for the involution $\psi$.}
\label{involution}
\end{figure}
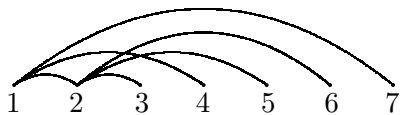

Here is an example for the involution $\psi$. Let $P$ be the linked partition given  in Figure \ref{pt}. Then  $\psi(P)$
is the linked partition in Figure \ref{involution}.

To conclude, we note that  Conjecture \ref{C-1} follows from Theorem
\ref{thm-f}. In fact,
Theorem \ref{thm-f} also implies Conjecture \ref{C-2}.
Let $P$ be  the linear representation of a linked partition of $[n]$, and let $S$ and $T$  be the sets of left-hand endpoints and right-hand endpoints of $P$, respectively.
It can be seen that $P$ is the front representation of a set partition of $[n]$
if and only if $S\cap T=\emptyset$.
Hence the set of front representations of  partitions of $[n]$
is the disjoint union of $L_n(S,T)$, where $(S,T)$ ranges over
pairs of disjoint subsets of $[n]$. 
For any two subsets $S$ and $T$ of $[n]$ with $S\cap T=\emptyset$,
if $L_n(S, T)$ is not empty, then we can apply Theorem \ref{thm-f} to deduce that  the crossing number and the nesting number have a symmetric joint distribution over $L_n(S,T)$.
Thus we have proved Conjecture \ref{C-2}.

\vskip 3mm \noindent {\bf Acknowledgments.}  This work was supported
by the 973 Project, the PCSIRT Project of
the Ministry of Education, the National Science Foundation of
China, the  Science Foundation of Hebei
Province and Grant No. BR2-231 of the Education Commission of
Hebei Province.

\end{document}